\documentclass[10pt,a4paper]{amsart}

\usepackage[all]{xy}
\usepackage[normalem]{ulem}
\usepackage{amsmath}
\usepackage{amsfonts}
\usepackage{amssymb}
\usepackage{amsthm}
\usepackage{graphicx}
\usepackage{microtype}
\usepackage[usenames,dvipsnames]{xcolor}
\usepackage{pgf,tikz}
\usepackage{mathrsfs}
\usetikzlibrary{arrows}
\usepackage{caption}
\usepackage{soul}
\usepackage[pdftex,plainpages=false,bookmarks=true,breaklinks=true,linkcolor=brown,citecolor=brown,colorlinks]{hyperref}
\usepackage{bm}
\usepackage{cancel} 
\usepackage{combelow}
\usepackage{here}

\newcommand{\prarrow}[2]{\ar@<0.5ex>[r]^-{#1} \ar@<-0.5ex>[r]_-{#2}}
\newcommand{\plarrow}[2]{\ar@<0.5ex>[l]^-{#1} \ar@<-0.5ex>[l]_-{#2}}
\newcommand{\pdarrow}[2]{\ar@<0.5ex>[d]^-{#1} \ar@<-0.5ex>[d]_-{#2}}
\newcommand{\puarrow}[2]{\ar@<0.5ex>[u]^-{#1} \ar@<-0.5ex>[u]_-{#2}}

\newtheorem{theorem}{Theorem}[section]
\newtheorem{corollary}[theorem]{Corollary}    
\newtheorem{lemma}[theorem]{Lemma}

\theoremstyle{definition}
   
\newtheorem{example}[theorem]{Example}    
\newtheorem{remark}[theorem]{Remark}

\author{Eiichi Matsuhashi}
\address{Department of Mathematics, Shimane University, Matsue, Shimane, 690-8504, Japan.}
\email{matsuhashi@riko.shimane-u.ac.jp}
\title{A non-$D$-continuum with weakly  infinite-dimensional closed set-aposyndetic Whitney levels}

\begin{document}

\begin{abstract} In this paper, we introduce  the new classes of continua;  {\it weakly infinite-dimensional closed set-aposyndetic continua}. With this notion, 
we show that there exists a  non-$D$-continuum  such that each positive Whitney level of the hyperspace  of the continuum is a weakly infinite-dimensional closed set-aposyndetic continuum. This result strengthens those of van Douwen and Goodykoontz \cite{douwen}, Illanes \cite{illanesapo}, and the main result of Illanes et al. \cite{moreon}.




\end{abstract}


\maketitle
\markboth{}{Eiichi Matsuhashi}

\section{Introduction}

In this paper, all spaces are metrizable and maps are continuous. Let $\mathbb{N}$ denote the set of all natural numbers, $\mathbb{Z}$ the set of all integers, and $\mathbb{R}$ the set of all real numbers. A compact metric space is called a $compactum$ and a $continuum$ means a connected compactum.  A continuum is said to be $decomposable$ if it is the union of two proper subcontinua of $X$.  If a continuum is not decomposable, it is said to be $indecomposable$. 

Let $X$ be
a compactum.  Then, $C(X)$ denotes the space of all subcontinua of $X$ with the topology generated by the Hausdorff metric. 
 $C(X)$ is called the $hyperspace$ $of$ $X$. 
A $Whitney$
$map$  is a map $\mu : C(X) \to [0,\mu(X)]$ satisfying  $\mu(\{x\})=0$ for each $x \in X$, and $\mu(A) < \mu(B)$  whenever $A, B \in C(X)$ and $A \subsetneq B$.  It is well-known  that for each Whitney map  $\mu: C(X) \to [0,\mu(X)]$   
and  each $t \in [0, \mu(X)]$,  $\mu^{-1}(t)$ is a continuum (\cite[Theorem 19.9]{illanes}). 

A topological property $P$, defined for continua, is: 

(a) a {\it Whitney property} if for any continuum $X$ possessing property $P$, every Whitney level of $C(X)$ also has property $P$; 

(b) a {\it Whitney reversible property} if whenever $X$ is a continuum such that every Whitney level has property $P$, it follows that $X$ itself has property $P$; 

(c) a {\it strong Whitney reversible property} if whenever $X$ is a continuum such that for some Whitney map $\mu$ for $C(X)$ and for every $t \in (0,1)$, $\mu^{-1}(t)$ has property $P$, then $X$ has property $P$; and 

(d)  a {\it sequential strong Whitney reversible property} provided
that whenever $X$ is a continuum such that there exist a Whitney map $\mu$
for $C(X)$ and a sequence $\{t_n\}_{n=1}^\infty$ in $(0,1)$ such
that $\lim\limits_{n\to\infty} t_n = 0$ and $ \mu^{-1}(t_n)$ has property
$P$ for each $n \ge 1$,  then $X$ has property $P$.

A comprehensive analysis of what was known on these subjects up to 1999 can be found in \cite[Chapter 8]{illanes}. 

Let $X$ be a compactum. Then, $X$ is said to be {\it weakly infinite-dimensional} if, for each countable  family $\{(A_i,B_i)\}_{i=1}^\infty$ of  pairs of disjoint closed subsets of $X$, there exists a family $\{L_i\}_{i=1}^\infty$ of closed subset of $X$ such that, for each $i \in \mathbb{N}$, $L_i$ is a partition in $X$ between $A_i$ and $B_i$. If $X$ is not weakly infinite dimensional, it is said to be {\it strongly infinite-dimensional}.  Each finite-dimensional space is weakly infinite-dimensional \cite[Proposition 6.1.2]{engelking} and the Hilbert cube is strongly infinite-dimensional \cite[Theorem 1.8.2]{engelking}.

Let $X$ be a continuum. Then,  $X$ is called a $D${\it -continuum} if for each pair of disjoint and nondegenerate subcontinua $A$ and $B$ of $X$, there exists a subcontinuum $L$ of $X$ such that $A\cap L\neq\emptyset\neq B\cap L$ and $(A\cup B)\setminus L\neq \emptyset$. In addition, if we require $A\setminus L\neq\emptyset\neq B\setminus L$ (resp., $B \setminus L \neq \emptyset$), then $X$ is called a $D^{*}${\it -continuum} (resp., a $D^{**}$-$continuum$). 

A continuum $X$ is called a {\it Wilder} continuum if for any three distinct points $x,y$ and $z$ of $X$, there exists a subcontinuum $L$ of $X$ such that $x\in L$ and $L$ contains exactly one of $y$ and $z$.  $X$ is said to be {\it continuum-wise Wilder} if
for each mutually disjoint subcontinua $A$, $B$ and $C$ of $X$, there exists a subcontinuum
$L$ of $X$ such that $A \cap L \neq \emptyset$ and exactly one of $B$ and
$C$ intersects $L$.
 $X$ is called a {\it closed set-wise Wilder continuum} if for each mutually disjoint closed subsets $A$, $B$ and $C$ of $X$, there exists a subcontinuum $L$ of $X$ such that $A \cap L \neq \emptyset$ and exactly one of $B$ and $C$ intersects $L$.

Let $X$ be a continuum and $n \in \mathbb{N}$. Then, $X$ is called an $n$-$aposyndetic\ continuum$ if for each $x \in X$ and each $n$-point subset $\{x_1,x_2,...,x_n\}$ of $X$ with  $\{x_1,x_2,...,x_n\} \subseteq X\setminus \{x\}$,  there exists a subcontinuum $T$ of $X$ such that $x \in {\rm Int}_X T \subseteq T \subseteq X \setminus \{x_1,x_2,...,x_n\}$, where ${\rm Int}_X T$ denotes the interior of $T$ in $X$. A $1$-aposyndetic continuum is called an $aposyndetic$ $continuum$.  If $X$ satisfies the condition that  for each $x \in X$ and each $0$-dimensional closed subset (resp., countable closed subset) $A$ of $X$ with $A \subseteq X\setminus \{x\}$,  there exists a subcontinuum $T$ of $X$ such that $x \in {\rm Int}_X T \subseteq T \subseteq X \setminus A$,  then $X$ is said to be {\it 0-dimensional closed set-aposyndetic (resp., countable close set-aposyndetic)}. Also, $X$ is called a $semi$-$aposyndetic \ continuum$ provided that for any two distinct points $x,y \in X$, there exists a subcontinuum $T$ of $X$ such that either $``x\in {\rm Int}_XT,~y\notin T"$ or $``x\notin T,~y\in {\rm Int}_XT."$

Note that every arcwise connected continuum is a Wilder continuum  and also a $D^*$-continuum. It is clear that every $D^*$-continuum is also a $D^{**}$-continuum, and every  $D^{**}$-continuum is a $D$-continuum. Furthermore, each Wilder continuum is a $D^{**}$-continuum (\cite[Proposition 2.2]{D**}). In  (\cite[Theorem~3.6]{D}), it is shown that  every semi-aposyndetic continuum is Wilder. Also, in the same paper, the authors showed that each 2-aposyndetic continuum is $D^*$ (\cite[Proposition 3.8]{D}). Since every nondegenerate indecomposable continuum has uncountably many composants (\cite[Theorem 11.15]{nadler1}), it follows  that every nondegenerate $D$-continuum is decomposable. For additional relationships between the classes of the above  continua and other classes of continua, refer to  \cite{camargo} and \cite{D}. Also, see \cite[Figure  1]{moreon} for further details. 

In this paper, we introduce  new classes of continua; {\it $n$-dimensional closed set-aposyndetic continua $(n \in \mathbb{N})$}, {\it  finite-dimensional closed set-aposyndetic continua},  and {\it weakly infinite-dimensional closed set-aposyndetic continua}. Let $X$ be a continuum and $n \in \mathbb{N}$. Then, $X$ is called an {\it n-dimensional closed set-aposyndetic continuum} if for each $x \in X$ and each $n$-dimensional closed subset $A$ of $X$ with $A \subseteq X\setminus \{x\}$,  there exists a subcontinuum $T$ of $X$ such that $x \in {\rm Int}_X T \subseteq T \subseteq X \setminus A$.  Also, $X$ is said to be {\it finite-dimensional closed set-aposyndetic (resp. weakly-infinite dimensional closed set-aposyndetic ) } if for each $x \in X$ and each finite-dimensional closed subset  (resp. weakly infinite-dimensional closed subset) $A$ of  $X$  with $A \subseteq X\setminus \{x\}$,  there exists a subcontinuum $T$ of $X$ such that $x \in {\rm Int}_X T \subseteq T \subseteq X \setminus A$. 


In \cite[p.47]{douwen}, Douwen and Goodykoontz provided an example of a non-semi-aposyndetic continuum  such that each positive Whitney level of its hyperspace is semi-aposyndetic. 
In \cite{illanesapo}, Illanes showed there exists a non-aposyndetic continuum $Z$ such that every positive Whitney level for $C(Z)$ is  finite-dimensional closed set-aposyndetic. In fact, in his example, every positive Whitney level for  $C(Z)$ is weakly infinite-dimensional closed set-aposyndetic. 
In \cite{moreon}, as the main result,   Illanes et al. constructed a non-$D$-continuum such that each positive Whitney level for the hyperspace of the continuum is both Wilder and $D^{*}$.  In Section 2, we provide one example that is stronger than all of these results. In Section 3,   we  provide an example of a finite-dimensional  closed set-aposyndetic continuum that is not weakly infinite-dimensional closed set-aposyndetic.  Also, we mention that  the example is useful for the construction of an $n$-dimensional closed set-aposyndetic continuum that is not 
$(n+1)$-dimensional closed set-aposyndetic for each $n \in \mathbb{N}$.






\section{Main result}

 In this section, by modifying the construction of the example given in \cite{illanesapo}, we show that there exists a non-$D$-continuum such that every positive Whitney level of the hyperspace of the continuum is a weakly infinite-dimensional closed set-aposyndetic continuum.
 
 If $(X,d)$ is a continuum, $A$ is a closed subset of $X$ and $\varepsilon > 0$,  then we denote the set $\{x \in X \ | \ {\rm there \ exists }\ a \in A \ {\rm such \ that } \ d(x,a) < \varepsilon \}$ by $N_d(A, \varepsilon)$. Also, $H_d$ denotes the Hausdorff metric on $C(X)$ induced  by $d$.
\begin{lemma}{\rm(\cite[Exercise 4.33 (b)]{nadler1})}
Let $(X,d)$ be a compactum and let $\mu : C(X) \to [0, \mu(X)]$ be a Whitney map. Then, for each $\varepsilon > 0$, there exists $\delta > 0$ such that if $A,B \in C(X)$ satisfy $B \subseteq N_d(A, \delta)$ and $|\mu(A) - \mu(B)| < \delta$, then $H_d(A,B) < \varepsilon$. 
\label{exe}
\end{lemma}
\begin{lemma}{\rm (\cite[Lemma 14.8.1]{nadler2})}
Let $X$ be a continuum, let $\mu:C(X) \to [0,\mu(X)]$ be a Whitney map and let $t \in (0, \mu(X))$. Then, for each distinct elements $A,B \in \mu^{-1}(t)$ with $A \cap B \neq \emptyset$,  there exists an arc $\mathcal{I}_{A,B} \subseteq \mu^{-1}(t)$ from $A$ to $B$ such that for each $L \in \mathcal{I}_{A,B}, ~ L \subseteq A \cup B$. Moreover, if $K$ is any given component of $A \cap B$, then   $\mathcal{I}_{A,B}$ may be chosen as above so that $K \subseteq L$ for each $L  \in \mathcal{I}_{A,B}$.
\label{lemmanadler}
\end{lemma}



A $dendrite$ is a locally connected continuum that contains no simple closed curve. A dendrite $D$ is called a $superdendrite$ if the set of all endpoints in $D$ is dense in $D$. A continuum  is said  to be {\it uniquely arcwise connected} provided that for every two of its points, there exists exactly one arc in the continuum joining these points. Since every dendrite is arcwise connected and does not contain a simple closed curve, it is easy to see the following result.

\begin{lemma}{\rm (Well-known)} Every dendrite is uniquely arcwise connected.
\label{uniquely}
\end{lemma}

Before Example \ref{main}, we introduce the following notation and terminology: If $X$ is a continuum and  $a \in X$, then $C(a,X)$ denotes the set $\{C \in C(X) : a \in C\}$. Also, $id_X$ denotes the identity map on $X$.  A continuum $Z$ is said to be $tree$-$like$ if for each $\varepsilon > 0,$ there exist a  tree 
 $T$ and a surjective map $f: Z \to T$ such that for each $t \in T$, the diameter of $f^{-1}(t)$ is smaller than $\varepsilon$.  Also, $Z$ is said to be $unicoherent$ if for each subcontinua $A, B \subseteq Z$ with $A \cup B=Z$, $A \cap B$ is connected. Furthermore, if each subcontinuum of $Z$ is unicoherent, then $Z$ is called a \textit{hereditarily unicoherent continuum}.

\begin{example}
There exists a non-$D$-continuum with weakly  infinite-dimensional closed set-aposyndetic Whitney levels. \label{main}
\end{example}

Let $d$ be the Euclidean   metric on $\mathbb{R}^3$.  Let $D$ be a superdendrite  contained in $[0,1]^2$ such that for each open subset  $O$ in $D$, $O$ contains a topological copy of the continuum $\bigcup_{n=1}^\infty\{(x,\frac{x}{n}) \in [0,1]^2 : 0 \le x \le \frac{1}{n}\}$. We may assume that there exist  two endpoints $a,b$ in $D$ such that $D \cap ([0,1] \times \{0\} )= \{a\}$ and $D \cap ([0,1] \times \{1\} )= \{b\}$. By Lemma \ref{uniquely}, there exists the unique subarc $I$ in $D$ from $a$ to $b$.  

For each $n \in \mathbb{N}$, let $L_n$ be the line segment in$[0,1]^2$ from $(1-\frac{1}{2^{(n-1)}},0)$ to  $(1-\frac{1}{2^{n}} ,1)$ if $n$ is odd , and the line segment in  $[0,1]^2$ from  $(1-\frac{1}{2^{(n-1)}},1)$ to  $(1-\frac{1}{2^{n}} ,0)$ if $n$ is even.  Next, for each $n  \in \mathbb{N}$, let  $D_n = \{(x,y,z) \in \mathbb{R}^3 : (x,z) \in L_n, \ (y,z) \in D \}$ and $D_{-n}= \{(x,y,z) \in \mathbb{R}^3: (-x,y,z) \in D_n \}$.  In addition, let $D_0=\{1\} \times D$ and $D_0' = \{-1\} \times D$.   Finally, define $Z=D_0 \cup D_0' \cup \bigcup _{n=1}^{\infty}(D_n \cup D_{-n} ) $. Note that $Z$ is a tree-like continuum. 

For each $n \in \mathbb{Z}$, let $a_n$, $b_n$ and $I_n$ be  points and subsets in $D_n$  corresponding to $a$,  $b$ and $I$ in $D$ respectively. Also, let $a_0$, $b_0$, and $I_0$ in $D_0$, and $a_0'$, $b_0'$, and $I_0'$ in $D_0'$ correspond to the points and subsets $a$, $b$, and $I$ in $D$, respectively. 

First, we show that $Z$ is not a $D$-continuum. Let   $L$ be a subcontinuum of $Z$ such that $I_0 \cap L \neq \emptyset \neq I_0' \cap L$. Then, it is easy to see that $\bigcup_{n \in \mathbb{Z}} I_n \subseteq L$. Hence, it follows that $I_0 \cup I_0' \subseteq L$. Therefore, we see that $Z$ is not a $D$-continuum.  

Let $t \in (0, \mu(Z))$.
We prove that $\mu^{-1}(t)$ is a weakly infinite-dimensional closed set-aposyndetic continuum. Let $S \in \mu^{-1}(t)$ and let $\mathcal{A}$ be a weakly infinite finite-dimensional closed subset of   $ \mu^{-1}(t) $ such that $S \notin \mathcal{A}$. 

\vspace{2mm}

{\bf Case 1.} $S \subsetneq D_0$.

Since $D_0$ is a dendrite, $S \notin \mathcal{A}$ and $\mathcal{A}$ is closed in $\mu^{-1}(t)$, by Lemma \ref{exe}, it is easy to see that   there exists a subcontinuum $\tilde{S}$ of $D_0$ such that $S \subseteq {\rm Int}_{D_0} \tilde{S}$ and $A \setminus \tilde S  \neq \emptyset$ for each $A \in \mathcal{A}$. Since each open subset of $D_0$ contains a topological copy of the continuum $\bigcup_{n=1}^\infty\{(x,\frac{x}{n}) \in [0,1]^2) : 0 \le x \le \frac{1}{n}\}$, it is not difficult to see that   $\mu^{-1}(t) \cap C(a_0, D_0)$ and $\mu^{-1}(t) \cap C(b_0, D_0)$ contain a topological copy of the Hilbert cube. Hence $\mu^{-1}(t) \cap C(a_0, D_0)$ and $\mu^{-1}(t) \cap C(b_0, D_0)$ are strongly infinite dimensional.  Therefore, $(\mu^{-1}(t) \cap C(a_0, D_0)) \setminus \mathcal{A} \neq \emptyset \neq (\mu^{-1}(t) \cap C(b_0, D_0)) \setminus \mathcal{A}$. Since $\mathcal{A}$ is weakly infinite-dimensional and $\mu^{-1}(t) \cap C(D_0)$ is homeomorphic to the Hilbert cube (see \cite[Theorem 4.8]{good}),  $(\mu^{-1}(t) \cap C(D_0)) \setminus \mathcal{A}$ is  connected (\cite{hadzii}. See also \cite[Problem 6.1.F]{engelking}). Hence, $(\mu^{-1}(t) \cap C(D_0)) \setminus \mathcal{A}$ is arcwise connected (\cite[Theorem 8.26]{nadler1}).  
 Therefore, we can see that there exist $\{F,G\} \subseteq \mu^{-1}(t) \cap C(D_0)$  and subarcs $\mathcal{J}_{S,F}$ and $\mathcal{J}_{S,G}$ of $\mu^{-1}(t) \cap C(D_0)$
such that 

\vspace{1mm}

(i)  $a_0 \in F$ and $b_0 \in G$, 

(ii) $\mathcal{J}_{S,F}$ is an arc from $S$ to $F$, and $\mathcal{J}_{S,G}$ is an arc from $S$ to $G$, and 

(iii) $\mathcal{A} \cap \mathcal{J}_{S,F}= \emptyset = \mathcal{A} \cap \mathcal{J}_{S,G}$.

\vspace{1mm}


Let $\pi : Z \to D_0$ be the natural projection. Note that $\pi|_{D_n}: D_n \to D_0$ is a homeomorphism for each $n \in \mathbb{N}$.  Let $\mathcal{S}=\{C \in \mu^{-1}(t) \cap C(D_0) : C \subseteq \tilde{S }\}$. 
 Also, for each $n \in \mathbb{N}$, let $\mathcal{S}_n= \{C \in \mu^{-1}(t) \cap C(D_0 \cup ( \bigcup_{i=n}^\infty D_i)) : C \subseteq \pi^{-1} (\tilde{S})  \} $. Note that $\mathcal{S}_n $ is  a neighborhood of $S$ in $\mu^{-1}(t)$ for each $ n \in \mathbb{N}$. 
Since $\mathcal{A}$ is closed in $\mu^{-1}(t)$ and $\mathcal{S} \cap \mathcal{A}= \emptyset$, we can easily see that   there exists $N_1 \in \mathbb{N}$ such that $\mathcal{S}_{N_1} \cap \mathcal{A}= \emptyset$.

Also,  imitating the argument in \cite[p.62]{illanesapo},  we can see the following:

\vspace{1mm}


$(\sharp)$ For each $n \in \mathbb{N}$, there exists a map $H_n: \mu^{-1}(t) \cap C(D_0) \to  \mu^{-1}(t) \cap C(D_n)$ such that 
$\{H_n\}_{n=1}^\infty$ uniformly converges to $id_{ \mu^{-1}(t) \cap C(D_0) }$.

\vspace{2mm}

Since $S \in {\rm Int}_{D_0} \tilde{S}$ and  $\{H_n\}_{n=1}^\infty $ uniformly converges to $id_{\mu^{-1}(t) \cap C(D_0)}$, there exists $N_2 \in \mathbb{N}$ such that if $n \in \mathbb{N}$ and  $n \ge N_2$, then $H_n(S) \in \mathcal{S}_n$.

For each $n \in \mathbb{N}$, let $\mathcal{J}^{H_n}_{S,F}=\{H_n(C) : C \in \mathcal{J}_{S,F}\}$ and $\mathcal{J}^{H_n}_{S,G}=\{H_n(C) : C \in \mathcal{J}_{S,G}\}$.  Note that $\mathcal{J}^{H_n}_{S,F}$ and $\mathcal{J}^{H_n}_{S,G}$ are subcontinua of $\mu^{-1}(t)$. 
Since  $ \mathcal{A}$ is closed in $\mu^{-1}(t)$,  $\mathcal{J}_{S,F} \cap \mathcal{A} = \emptyset = \mathcal{J}_{S,G} \cap \mathcal{A}$, and $\{H_n\}_{n=1}^\infty $ uniformly converges to $id_{\mu^{-1}(t) \cap C(D_0)}$,  there exists $N_3 \in \mathbb{N}$ such that if   $n \in\mathbb{N}$ and   $n \ge  N_3$, then $\mathcal{J}^{H_n}_{S,F} \cap \mathcal{A} = \emptyset = \mathcal{J}^{H_n}_{S,G} \cap \mathcal{A}$.

\vspace{2mm}

\textbf{Claim.} If $\{A_n\}_{n=1}^\infty$,  $\{B_n\}_{n=1}^\infty \subseteq C(D_0)$,  $\lim A_n  = F$,  and $\lim B_n =G$, then there exists $ M \in \mathbb{N}$ such that $A_m \cap A_n \neq \emptyset \neq B_m \cap B_n $ for each $m,n > M$.   

Proof of  Claim.  We only prove that there exists $M \in \mathbb{N}$ such that $A_m \cap A_n \neq \emptyset$ for each $m,n > M$. If not, we may assume $A_{2n-1} \cap A_{2n} = \emptyset$ for each $n \in \mathbb{N}$. Since $D_0$ is a dendrite, there exists $\delta_0 > 0$ such that if $x,y \in D_0$ and $d(x,y) < \delta_0$, then there exists a  subarc $J$ of $D_0$ from $x$ to $y$ such that ${\rm diam} J < \frac{{\rm diam} F}{8}$ (see \cite[Exercise 8.43]{nadler1}).  Since $\lim A_n = F$, there exists $M' \in \mathbb{N}$ such that  ${\rm diam}A_n > \frac{7}{8} {\rm diam }F$  and $H_d(A_{2n-1}, A_{2n}) <  \delta_0 $ for each $n  \ge M'$, where $H_d$ denotes the Hausdorff metric on $C(Z)$.  
Then, it is easy to see that there exist  pairwise disjoint subarcs $J_0$ and $J_1$ of $D_0$ such that $J_0 \cap A_{2M'-1} \neq \emptyset \neq J_0 \cap  A_{2M'} $ and $J_1 \cap A_{2M'-1} \neq \emptyset \neq J_1 \cap  A_{2M'}$. 
 Then, we  see that $A_{2M'-1} \cup A_{2M'} \cup J_0 \cup J_1$ contains a simple closed curve. This is a contradiction since $D_0$ is a dendrite. This completes the proof of Claim.

\vspace{2mm}

By Order Arc Theorem (see \cite[Exercise 5.25]{nadler1}), for each $n \in \mathbb{N}$ , we can take $F_n, G_n \in \mu^{-1}(t) \cap C(D_n)$  such that $a_n \in F_n \subseteq   \pi_n^{-1}(F)$ or $a_n  \in   \pi_n^{-1}(F) \subseteq F_n$, and $b_n \in G_n \subseteq   \pi_n^{-1}(G)$ or $b_n  \in   \pi_n^{-1}(G) \subseteq G_n$. By Lemma \ref{exe} and $(\sharp)$, we see that $\lim F_n = F =  \lim H_n(F)$ and  $\lim G_n = G =  \lim H_n(G)$. Hence, using Claim, we see that there exists $N_4 \in \mathbb{N}$ such that  $F_n \cap H_n(F)  \neq \emptyset \neq G_n \cap  H_n (G) $ for each $n > N_4$.

Let $\varepsilon >0$.  For this $\varepsilon$, take $\delta > 0$ as in the statement of Lemma \ref{exe}. We may assume that $\delta < \varepsilon$. Take  $N_5 \in \mathbb{N}$ such that $N_5 > N_4$  and if $n \in \mathbb{N}$ and $n \ge N_5$, then $H_d(F_n \cup H_n(F), F) < \delta$, $H_d(G_n \cup H_n(G), G) < \delta$, $t \le \mu (F_n \cup H_n(F))  \le t+ \delta $ and $t \le \mu (G_n \cup H_n(G))  \le  t+ \delta $.
 By  Lemma \ref{lemmanadler}, there exist subarcs $\mathcal{J}_{F_n, H_n(F)}$ of $\mu^{-1}(t) \cap C (F_n \cup H_n(F))$ from $F_n$ to $H_n(F)$ and $\mathcal{J}_{G_n, H_n(G)}$ of $\mu^{-1}(t) \cap C (G_n \cup H_n(G))$ from $G_n$ to $H_n(G)$ (or $\mathcal{J}_{F_n, H_n(F)}$ and $\mathcal{J}_{G_n, H_n(G)}$ might be one point sets).  Then, for each $C \in \mathcal{J}_{F_n, H_n(F)}$, $H_d(F,C) <  \varepsilon$, and for each $C \in \mathcal{J}_{G_n,H_n(G)}$, $H_d(G,C) <  \varepsilon$.  Hence, we see that $\lim \mathcal{J}_{F_n, H_n(F)} = \{F\}$ and  $\lim \mathcal{J}_{G_n, H_n(G)} = \{G\}$ (this limit is with respect to the metric $H_{H_d}$).   In particular, we see that  $\lim F_n = F$ and $\lim G_n =G$. Hence,  $\lim (F_n \cup F_{n+1})= F$ and $\lim (G_n \cup G_{n+1}) =G$. Also, note that there exists $N_6 \in \mathbb{N}$ such that  if $n \in \mathbb{N}$ and $n \ge N_6$, then $\mathcal{J}_{F_n, H_n(F)} \cap \mathcal{A} = \emptyset = \mathcal{J}_{G_n, H_n(G)} \cap \mathcal{A}$.




For each $n \in \mathbb{N}$, let $\mathcal{K}_{F_n} =  \mu^{-1}(t) \cap C(F_n \cup F_{n+1}) $   and  $\mathcal{K}_{G_n} =  \mu^{-1}(t) \cup C(G_n \cap G_{n+1})  $.   By Lemma \ref{exe} and the fact that $\lim (F_n \cup F_{n+1})= F$ and $\lim (G_n \cup G_{n+1}) =G$, we can easily see that $\lim \mathcal{K}_{F_n}= \{F\}$  and $\lim \mathcal{K}_{G_n}= \{G\}$. Since   $\lim \mathcal{K}_{F_n}= \{F\}$, $\lim \mathcal{K}_{G_n}= \{G\}$ and $F,G \notin \mathcal{A}$, there exists $N_7 \in \mathbb{N}$ such that 
 if $n \in \mathbb{N}$ and $n \ge N_7$, then $\mathcal{K}_{F_n} \cap \mathcal{A} = \emptyset = \mathcal{K}_{G_n} \cap \mathcal{A}$. 


 Let $N=\max\{N_1,N_2,N_3, N_5,N_6, N_7\}$, and let $\mathcal{T}=  \mathcal{J}_{S,F} \cup \mathcal{J}_{S,G} \cup \mathcal{S}_N \cup (\bigcup_{n=N}^{\infty}( \mathcal{J}^{H_n}_{S,F} \cup \mathcal{J}^{H_n}_{S,G}\cup  \mathcal{J}_{F_n, H_n(F)} \cup \mathcal{J}_{G_n, H_n(G)} \cup \mathcal{K}_{F_n} \cup \mathcal{K}_{G_n}) )$.  Note that $\mathcal{T} \cap \mathcal{A}= \emptyset$.

Since ${\rm Cl}_{\mu^{-1}(t)} (\mathcal{S}_N \setminus \mathcal{S})  = \mathcal{S}_N$, $\lim \mathcal{J}^{H_n}_{S,F}= \mathcal{J}_{S,F}$,  $\lim \mathcal{J}^{H_n}_{S,G}= \mathcal{J}_{S,G}$,  \ $\lim \mathcal{J}_{F_n, H_n(F)}= \{F\}$, $\lim \mathcal{J}_{G_n, H_n(G)}= \{G\}$, $\lim \mathcal{K}_{F_n} =\{F\}$ and  $\lim \mathcal{K}_{G_n} =\{G\}$, ${\rm Cl}_{\mu^{-1}(t)}((\mathcal{S}_N \setminus \mathcal{S}) \cup (\bigcup_{n=N}^{\infty}(    \mathcal{J}^{H_n}_{S,F} \cup \mathcal{J}^{H_n}_{S,G}\cup  \mathcal{J}_{F_n, H_n(F)} \cup \mathcal{J}_{G_n, H_n(G)} \cup \mathcal{K}_{F_n} \cup \mathcal{K}_{G_n})))=\mathcal{T}$. Note that $(\mathcal{S}_N \setminus \mathcal{S}) \cup (\bigcup_{n=N}^{\infty}(    \mathcal{J}^{H_n}_{S,F} \cup \mathcal{J}^{H_n}_{S,G}\cup  \mathcal{J}_{F_n, H_n(F)} \cup \mathcal{J}_{G_n, H_n(G)} \cup \mathcal{K}_{F_n} \cup \mathcal{K}_{G_n}))$ is connected. Hence, $\mathcal{T}$ is a continuum.  Also, it is easy to see that  $S \in {\rm Int}_{\mu^{-1}(t)} \mathcal{S} _N     \subseteq  {\rm Int}_{\mu^{-1}(t)} \mathcal{T} \subseteq \mathcal{T} \subseteq \mu^{-1}(t) \setminus  \mathcal{A} $.   Hence,   $\mathcal{T}$ is a subcontinuum  of $\mu^{-1}(t)$ such that $S \in \rm{Int}_{\mu^{-1}(t)}\mathcal{T} \subseteq \mathcal{T} \subseteq \mu^{-1}(t) \setminus \mathcal{A}$.

\vspace{2mm}

{\bf Case 2. } $S \subsetneq D_0'$. 

In this case, as in the argument of Case 1, 
 we can find a subcontinuum $\mathcal{T}$ of $\mu^{-1}(t)$ such that $S \in \rm{Int}_{\mu^{-1}(t)}\mathcal{T} \subseteq \mathcal{T} \subseteq \mu^{-1}(t) \setminus \mathcal{A}$.

\vspace{2mm} 

{\bf Case 3. } $S \cap D_0 \neq \emptyset \neq S \cap D_0'$, or $S \cap D_0=\emptyset$ and $S \cap D_0' \neq \emptyset \neq S \cap (Z \setminus (D_0 \cup D_0'))$, or $S \cap D'_0=\emptyset$ and $S \cap D_0 \neq \emptyset \neq S \cap (Z \setminus (D_0 \cup D_0'))$, or $S \subseteq Z \setminus (D_0 \cup D_0') $, or $S=D_0$, or $S=D_0'$ .

In Case 3, the condition that $\mathcal{A}$ is weakly infinite-dimensional is not utilized. 
We only prove the case where $S \cap D_0 \neq \emptyset \neq S \cap D_0'$ (the proofs for the other cases can be obtained in a similar or a simpler way). 
 Let $S_0= S \cap D_0$ and $S_0'= S \cap D_0'$.  Since $Z$ is tree-like, $Z$ is hereditarily unicoherent (see \cite[p.232]{nadler1}). Hence, $S_0$ and $S_0'$ are continua.  Note that   $I_0 \subseteq S_0$ and   $I_0' \subseteq S_0'.$  Since $\mathcal{A}$ is closed in $\mu^{-1}(t)$ and $S \notin \mathcal{A}$, by applying Lemma \ref{exe}, we can take   subcontinua $\tilde{S_0}$ of $D_0$,  $\tilde{S_0'}$ of $D_0'$,  $N' \in \mathbb{N} $ and a  connected closed neighborhood $V$ of  $S \cap (Z \setminus (D_0 \cup D_0' ))$ in $ Z \setminus (D_0 \cup D_0')$  such that $\mu(\tilde{S_0}), \mu(\tilde{S'_0}) \in (0,t)$, $S_0 \subseteq {\rm Int}_{D_0} \tilde{S_0}$, $S_0 \subseteq {\rm Int}_{D_0'} \tilde{S_0'}$, 
 $ ({\rm Cl}_Z V) \cap D_0  \neq \emptyset  \neq ({\rm Cl}_Z V) \cap D_0'   $, $  ({\rm Cl}_Z V) \cap D_0  \subseteq S_0 $, $ ({\rm Cl}_Z V) \cap D_0'  \subseteq S_0' $,   and $A \setminus (V \cup (\tilde{S_0} \times [1-\frac{1}{2^{N'}},1]) \cup (\tilde{S_0} \times [-1 , -(1-\frac{1}{2^{N'}})]))\neq \emptyset$ for each $A \in \mathcal{A}$. Then, $\mathcal{T}=\{C \in \mu^{-1}(t) : C \subseteq V \cup (\tilde{S_0} \times [1-\frac{1}{2^{N'}},1]) \cup (\tilde{S_0'} \times [-1 , -(1-\frac{1}{2^{N'}})])\}$ is a subcontinuum of $\mu^{-1}(t)$ such that $S \in \rm{Int}_{\mu^{-1}(t)}\mathcal{T} \subseteq \mathcal{T} \subseteq \mu^{-1}(t) \setminus \mathcal{A}$.  

\vspace{2mm}

In every case, we can find a subcontinuum of $\mathcal{T}$ of  $\mu^{-1}(t)$ such that $S \in \rm{Int}_{\mu^{-1}(t)}\mathcal{T} \subseteq \mathcal{T} \subseteq \mu^{-1}(t) \setminus \mathcal{A}$. Hence, $\mu^{-1}(t)$ is weakly infinite-dimensional closed set-aposyndetic.

\begin{corollary} {\rm (cf. \cite[p.47]{douwen}, \cite[EXAMPLE]{illanesapo} and \cite[Example 3.2]{illanes})} Each of the following properties is not a Whitney reversible property:  being $D$, being $D^*$, being $D^{**}$, being Wilder, being  semi-aposyndetic, being countable closed set-aposyndetic, being $n$-aposyndetic $(n \in \mathbb{N} \cup \{0\})$,   being $n$-dimensional closed set-aposyndetic $(n \in \mathbb{N} \cup \{0\})$, being finite-dimensional closed set-aposyndetic, and being weakly infinite-dimensional closed set-aposyndetic.   
\end{corollary}

\begin{remark}
Being closed set-wise Wilder and being continuum-wise Wilder are  sequential strong Whitney reversible properties (see \cite[Corollaries 4.3 and 5.2]{illanesapo}).

\end{remark}

 \section{Examples}
 
We first provide an example of a finite-dimensional  closed set-aposyndetic continuum that is not weakly infinite-dimensional closed set-aposyndetic. 

\begin{example}

Let $X$ be the Hilbert cube.  Consider a family $\{X_n\}_{n=1}^{\infty}$ of subsets of $X$ such that:

(1) for each $n \in \mathbb{N}$, $X_n$ is homeomorphic to $[0,1]^n$,

(2) for each $n,m \in \mathbb{N}$, $X_n \cap X_m \neq \emptyset$ if and only if $|n - m| \leq 1$,

(3) for each $n \in \mathbb{N}$, there exists a point $p_n \in X$ such that $X_n \cap X_{n+1} = \{p_n\}$, and

(4) there exists a point $p \in X \setminus \bigcup_{n=1}^\infty X_n$ such that ${\rm Cl}_{X}(\bigcup_{n=1}^\infty X_n) \setminus \bigcup_{n=1}^\infty X_n = \{p\}$.

Take topological copies $F,G \subseteq X$ of  $\{ p\} \cup (\bigcup_{n=1}^\infty X_n)$ such that $F \cap G = \emptyset$. Note that $F$ and $G$ are weakly infinite-dimensional (see \cite[Theorem 6.1.11]{engelking}).  


Let $Y_0 = X \times \{0\}$.  
For each $n \in \mathbb{N}$, let  $g_n: X \to [\frac{1}{n+1},\frac{1}{n}]$ be a map such that $g_n^{-1}(\frac{1}{n})= F$ and   $g_n^{-1}(\frac{1}{n+1})= G$ if $n$ is odd, and $g_n^{-1}(\frac{1}{n})= G$ and   $g_n^{-1}(\frac{1}{n+1})= F$ if $n$ is even.  Also, for each $n \in \mathbb{N}$, let $Y_n=\{((x,y) \in X \times [0,1]) : y=f_n(x)\}$. Finally, let  $Z= \bigcup_{n=0}^\infty Y_n$. By a similar argument in Example \ref{main}, we can see that $Z$ is a finite-dimensional closed set-aposyndetic continuum  that is not weakly infinite-dimensional closed set-aposyndetic.

\end{example}

\begin{remark}
Let  $n \in \mathbb{N}$. In the previous example, by replacing $F$ and $G$ with a topological copy of $[0,1]^{n+1}$, we can obtain an $n$-dimensional closed set-aposyndetic continuum that is not $(n+1)$-dimensional closed set-aposyndetic. Similarly,   by replacing $F$ and $G$ with a Cantor set,  we can obtain a countable  closed set-aposyndetic continuum that is not 0-dimensional closed set-aposyndetic.  Furthermore, by replacing $F$ and $G$ with an  $(n+1)$-point set,  we can obtain an $n$-aposyndetic continuum that is not $(n+1)$-aposyndetic.  
\end{remark}

\section{Acknowledgements}

The  author was partially supported by JSPS KAKENHI Grant number 21K03249.

\end{document}